\documentclass[10pt]{article}
\pdftrailerid{}
\usepackage{subcaption}
\usepackage{comment}

\usepackage[utf8]{inputenc}
\usepackage[T1]{fontenc}
\usepackage{amsmath,amssymb,dsfont}
\numberwithin{equation}{section}
\usepackage{microtype}
\usepackage{graphicx,tikz,pgfplots}
\usetikzlibrary{positioning}
\graphicspath{{images/}}
\pgfplotsset{compat=newest}
\usepackage[hyperref,amsmath,thmmarks]{ntheorem}
\usepackage{bm}
\usepackage{aliascnt}
\usepackage{float}
\usepackage[a4paper,centering,bindingoffset=0cm,marginpar=2cm,margin=2.5cm]{geometry}
\usepackage[pagestyles]{titlesec}
\usepackage[font=footnotesize,format=plain,labelfont=sc,textfont=sl,width=0.75\textwidth,labelsep=period]{caption}
\usepackage{mathtools}

\makeatletter
\begingroup
\endlinechar=-1\everyeof{\noexpand}
\edef\x{\endgroup\def\noexpand\bibinputs{\@@input|"kpsewhich infmath.bib" }}\x
\makeatother
\ifx\bibinputs\empty\else\def\havecscbib{}\fi

\usepackage[backend=biber,maxnames=10,backref=true,hyperref=true,giveninits=true,safeinputenc]{biblatex}
\ifdefined\havecscbib
    \ifdefined\extract
    \else
    \fi
\else
    \ifdefined\extract
      \errmessage{cannot extract without detected csc bibliography}
    \fi
\fi

\DefineBibliographyStrings{english}{%
	backrefpage = {cited on page},
	backrefpages = {cited on pages},
}

\DeclareFieldFormat[report]{title}{``#1''}
\DeclareFieldFormat[book]{title}{``#1''}
\AtEveryBibitem{\clearfield{url}}
\AtEveryBibitem{\clearfield{note}}

\titleformat{\section}[block]{\large\sc\filcenter}{\thesection.}{0.5ex}{}[]
\titleformat{\subsection}[runin]{\bf}{\thesubsection.}{0.5ex}{}[.]

\usepackage[pdftex,colorlinks=true,linkcolor=blue,citecolor=green,urlcolor=blue,bookmarks=true,bookmarksnumbered=true]{hyperref}
\hypersetup
{
    bookmarksnumbered,
    breaklinks,
    pdfborderstyle=,
    colorlinks,
    citecolor=[rgb]{0,.65,0},
    linkcolor=[rgb]{0,0,.5},
    urlcolor=[rgb]{0,0,.5},
    pdfauthor={Authors},
    pdfsubject={Subject},
    pdftitle={Title},
    pdfkeywords={Keywords}
}

\newpagestyle{headers}
{
	\headrule
	\sethead[\footnotesize\thepage][\footnotesize\sc Authors][]{}{\footnotesize\sc Title}{\footnotesize\thepage}
	\setfoot{}{}{}
}
\newaliascnt{proposition}{lemma}

\aliascntresetthe{proposition}

\newaliascnt{corollary}{lemma}

\aliascntresetthe{corollary}

\newaliascnt{theorem}{lemma}

\aliascntresetthe{theorem}

\theorembodyfont{\normalfont}
\newaliascnt{definition}{lemma}

\aliascntresetthe{definition}

\newaliascnt{assumption}{lemma}

\aliascntresetthe{assumption}

\newaliascnt{example}{lemma}
\newtheorem{example}[example]{Example}
\aliascntresetthe{example}

\theoremstyle{nonumberplain}
\theoremseparator{:}
\theoremheaderfont{\normalfont\itshape}

\newtheorem{remark}{Remark}

\theoremsymbol{\ensuremath{\square}}

\newcommand{\R}{\mathds{R}}

\let\RE\Re
\let\Re=\undefined
\DeclareMathOperator{\Re}{\RE e}
\let\IM\Im
\let\Im=\undefined
\DeclareMathOperator{\Im}{\IM m}

\newcommand{\abs}[1]{\left|#1\right|}
\newcommand{\norm}[1]{\left\|#1\right\|}
\newcommand{\set}[1]{\left\{#1\right\}}

\newcommand{\e}{\mathrm e}
\let\ii\i
\renewcommand{\i}{\mathrm i}

\ifdefined\final
	\newcommand{\commentK}[1]{}
	\newcommand{\commentS}[1]{}
	\newcommand{\commentJ}[1]{}
	\newcommand{\commentO}[1]{}
	\newcommand{\commentC}[1]{}
	\newcommand{\CK}[1]{}
\else
	\newcommand{\commentK}[1]{{\color{red}{\bf Kemal:} #1}}
	\newcommand{\commentS}[1]{{\color{orange}{\bf Christina:} #1}}
	\newcommand{\commentJ}[1]{{\color{green}{\bf Jikai:} #1}}
	\newcommand{\commentO}[1]{{\color{blue}{\bf Otmar:} #1}}
	\newcommand{\commentC}[1]{{\color{brown}{\bf Clemens:} #1}}
	\newcommand{\CK}[1]{{\color{brown}#1}}
\fi

\newcommand{\vx}{{\bf x}}

\def\vek#1{\bm{#1}}

\title{Using the Transport of Intensity and the Transport of Phase Equation for Phase Retrieval}
\author{Clemens Kirisits$^1$\\{\footnotesize\href{mailto:clemens.kirisits@univie.ac.at}{clemens.kirisits@univie.ac.at}}
    \and Kemal Raik$^1$\\{\footnotesize\href{mailto:kemal.raik@univie.ac.at}{kemal.raik@univie.ac.at}}
	\and Otmar Scherzer$^{1,2,3}$\\{\footnotesize\href{mailto:otmar.scherzer@univie.ac.at}{otmar.scherzer@univie.ac.at}}
	\and Christina Strohmenger$^{1}$\\{\footnotesize\href{mailto:christina.ulrike.strohmenger@univie.ac.at}{christina.ulrike.strohmenger@univie.ac.at}}
	\and Jikai Yan$^{1}$\\{\footnotesize\href{mailto:jikai.yan@univie.ac.at}{jikai.yan@univie.ac.at}}
}
\date{}

\begin{document}

\maketitle
\thispagestyle{empty}
\begin{center}
	\hspace*{5em}
	\parbox[t]{12em}{\footnotesize
		\hspace*{-1ex}$^1$Faculty of Mathematics\\
		University of Vienna\\
		Oskar-Morgenstern-Platz 1\\
		A-1090 Vienna, Austria}
	\hfil
	\parbox[t]{17em}{\footnotesize
		\hspace*{-1ex}$^2$Johann Radon Institute for Computational\\
		\hspace*{1em}and Applied Mathematics (RICAM)\\
		Altenbergerstraße 69\\
		A-4040 Linz, Austria}
\end{center}

\begin{center}
\parbox[t]{19em}{\footnotesize
\hspace*{-1ex}$^3$Christian Doppler Laboratory\\
for Mathematical Modeling and Simulation\\
of Next Generations of Ultrasound Devices (MaMSi)\\
Oskar-Morgenstern-Platz 1\\
A-1090 Vienna, Austria}
\end{center}
\begin{abstract}
We investigate the transport of intensity equation (TIE) and the transport of phase equation (TPE) for solving the phase retrieval problem. Both the TIE and the TPE are derived from the paraxial Helmholtz equation and relate phase information to the intensity. The TIE is usually favored since the TPE is nonlinear. The main contribution of this paper is a discussion of preferential use of either one of the two equations or potential benefits of a hybrid use. Moreover, we discuss the solution of the TPE with the method of characteristics and with viscosity methods. Both the TIE and the viscosity method are numerically implemented with finite element methods.
\end{abstract}
\section{Introduction}
\label{sec:introduction}
The phase retrieval problem consists in computing the phase $\phi:\mathbb{R}^3\to\mathbb{R}$ of a complex valued function 
\begin{equation} \label{eq:AI}
    A:\mathbb{R}^3 \to \mathbb{C}, \quad A = \sqrt{I} \exp(\i \phi),
\end{equation}
given measurements of the intensity $I:\mathbb{R}^3\to\mathbb{R}$. Since typically the phase is more difficult to measure than the intensity but the phase describes the specimen more appropriately, the phase retrieval problem has several applications in the fields of microscopy \cite{Tutorial}, holography and crystallography (see for instance \cite{GeiDui71}).
There exist experimental setups, such as interferometers \cite{Fit00, DavNohSolZie02, GeiDui71} and computational algorithms for phase retrieval, examples for the latter being the Gerchberg-Saxton algorithm \cite{GerSax72} and the Hybrid Input-Out method proposed by Fienup \cite{Fie78}. 

One approach for solving the phase retrieval problem in the case that the illuminating wave is traveling primarily in the $z$-direction consists in solving the \textbf{Transport of Intensity Equation} (TIE)
\begin{equation} \label{eq:TIE}
\nabla_{\vx}\cdot\left[I(\vx,z)\nabla_{\vx}\phi(\vx,z)\right]=kI_z(\vx,z)
\end{equation}
for the phase $\phi$. Throughout this paper we let $\vx = (x,y) \in \R^2$. The equation was originally proposed by Teague \cite{Tea83} and has been successfully applied in microscopy to reconstruct the phase $\phi$ from intensity measurements $I$ (see \cite{BosFroNilSagUns16}). Several numerical methods to solve this equation have been studied in the literature: In \cite{Tea83} an auxiliary function $\nabla_{\vx}\psi=I\nabla_{\vx}\phi$ is used to convert the problem into finding the solution of a Poisson equation 
\begin{equation} \label{eq:poissontrick}
   \Delta_{\vx}\psi= kI_z.
\end{equation} 
The solution can then be plugged into a second Poisson equation,  $\Delta_{\vx}\phi=\nabla_{\vx}\cdot\left[\frac{1}{I}\nabla_{\vx}\psi \right]$, which may be solved with a Green's function. The implementation is based on FFT methods to compute the inverse of the Laplacian operator. FFT methods inherently assume periodic boundary conditions. \cite{GurRobNug95} solves the TIE under Neumann boundary conditions making use of a Zernike polynomial decomposition. Bostan et al.\ \cite{BosFroNilSagUns16} use regularisation methods to find the solution of the TIE under periodic boundary conditions. They also rewrite the TIE as a Poisson equation. As we show below, boundary conditions are important to determine the correct phase (see \autoref{subsec:ex1} and \autoref{subsec:ex2}). One alternative approach, which avoids specification of boundary conditions is by considering the TIE in free space with solution in a Sobolev space. A difficulty associated with this strategy is that $I$ has compact support, such that ellipticity of the differential operator in \autoref{eq:TIE} is violated. 
Moreover, Sobolev spaces typically require decay at infinity, which means that the phase has to be zero at infinite, which seems an unrealistic assumptions. See \cite{Eng97}, for the analysis of the decay rate of the solution of the integral equation with the associated Green function as kernel. We emphasize that \autoref{eq:poissontrick} of \cite{Tea83} is an attempt to avoid using the Green's function with diffusion coefficient $I$. 

The starting point for this paper is the following observation: In \autoref{sec:derivation} we derive the TIE from the Helmholtz equation assuming that the wave travels mainly in the $z$-direction. In doing so, a second equation appears, the {\bf Transport of Phase Equation} (TPE)
\begin{equation} \label{eq:TPE}
2k\phi_z-\|\nabla_{\vx}\phi\|^2=-\frac{\Delta_{\vx}\sqrt{I}}{\sqrt{I}} =:-\hat{I}.
\end{equation}

Both the TIE and the TPE relate the phase to the intensity of the wave. In this paper we will discuss the question of whether the solutions differ if we use one or the other equation for the reconstruction.
Furthermore, in applications such as probes being analyzed under the microscope, the domain of interest 
is bounded.  Therefore, both equations require appropriate boundary data in order to be uniquely solvable. In addition to the above mentioned settings, in \cite{Petrucelli} homogeneous Dirichlet boundary conditions are used for the reconstructions; \cite{Rod90} considers the TIE with Neumann boundary conditions.

In the following, we give an overview on the information that is required in order to reconstruct the phase based on the TIE and the TPE, respectively. For the simplicity of presentation we consider the equations on the three dimensional cube 
$$\Omega= [0,1]^3 \text{ and } \Omega_z = [0,1]^2 \times \lbrace z\rbrace.
$$ 
Lastly, we write $\Gamma_z:=\partial \Omega_z$ for the boundary of the two-dimensional square $\Omega_z$. We assume that $I$ is known in all of $\Omega.$

We distinguish two reconstruction situations, which are illustrated in \autoref{fig:introboundary}:
\begin{description}
    \item [Case 1] Given boundary data of $\phi$ on $\Gamma_{z^*}$, we can reconstruct $\phi$ in $\Omega_{z^*}$ by solving the TIE. If phase information is available on $\Gamma_z$ for every $z \in [0,1]$, then we can reconstruct $\phi$ on the whole domain $\Omega$. 
    \item [Case 2] If $\phi$ is given on $\Omega_{z^*}$, then it can be reconstructed at least in a neighborhood of $\Omega_{z^*}$ by solving the TPE with the method of characteristics (see \autoref{sec:TPE}). In particular, given $\phi$ on $\Gamma_{z^*}$, we can compute $\phi$ in $\Omega_{z^*}$ with the TIE and then reconstruct $\phi$ locally around $\Omega_{z^*}$ with the TPE. 
\end{description}

The hybrid constellation of the TIE and the TPE has the advantage that only a single set of boundary data of the phase (see \autoref{subsection: charac} for details) are required to retrieve $\phi$ at least locally. On the other hand, reconstructing $\phi$ in $\Omega$ using the TIE alone would require boundary data on each boundary $\Gamma_z$, $z \in [0,1]$ as pointed out in Case 1. 
However, if the TPE is solved computationally based on the viscosity method (see \autoref{subsection: visco}), we need again boundary data on $\Gamma_z$ for $z\in[0,1]$, that is, the same boundary data as for the pure TIE approach. 

\tikzset{every picture/.style={line width=0.75pt}}        
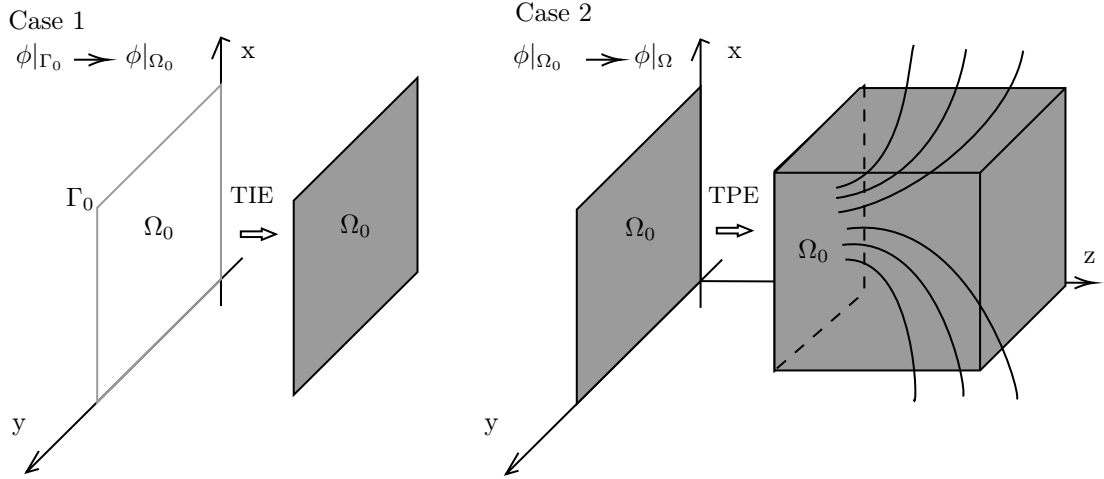
\begin{figure}[h]
    \centering
     
\begin{tikzpicture}[x=0.75pt,y=0.75pt,yscale=-1,xscale=1,scale=0.9]

	\draw  (130.54,173.69) -- (130.54,23.69)(22.54,266.69) -- (142.54,146.69) (125.54,35.69) -- (130.54,23.69) -- (135.54,25.69) (29.54,264.69) -- (22.54,266.69) -- (29.54,254.69)  ;

	\draw    (386.54,159.69) -- (612.92,160.68) ;
	\draw [shift={(614.92,160.69)}, rotate = 180.25] [color={rgb, 255:red, 0; green, 0; blue, 0 }  ][line width=0.75]    (10.93,-3.29) .. controls (6.95,-1.4) and (3.31,-0.3) .. (0,0) .. controls (3.31,0.3) and (6.95,1.4) .. (10.93,3.29)   ;

	\draw  [color={rgb, 255:red, 155; green, 155; blue, 155 }  ,draw opacity=1 ] (130.54,158.69) -- (61.88,227.35) -- (61.88,118.73) -- (130.54,50.07) -- cycle ;

	\draw  (396.54,174.69) -- (396.54,24.69)(288.54,267.69) -- (408.54,147.69) (391.54,36.69) -- (396.54,24.69) -- (401.54,26.69) (295.54,265.69) -- (288.54,267.69) -- (295.54,255.69)  ;

	\draw  [fill={rgb, 255:red, 155; green, 155; blue, 155 }  ,fill opacity=1 ] (239.54,154.69) -- (170.88,223.35) -- (170.88,114.73) -- (239.54,46.07) -- cycle ;

	\draw  [fill={rgb, 255:red, 155; green, 155; blue, 155 }  ,fill opacity=1 ] (396.54,159.69) -- (327.88,228.35) -- (327.88,119.73) -- (396.54,51.07) -- cycle ;

	\draw   (141.54,131.94) -- (153.17,131.94) -- (153.17,130.36) -- (160.92,133.53) -- (153.17,136.69) -- (153.17,135.11) -- (141.54,135.11) -- cycle ;

	\draw   (405.54,129.94) -- (417.17,129.94) -- (417.17,128.36) -- (424.92,131.53) -- (417.17,134.69) -- (417.17,133.11) -- (405.54,133.11) -- cycle ;

	\draw    (49,36) -- (66,36) ;
	\draw [shift={(68,36)}, rotate = 179.4] [color={rgb, 255:red, 0; green, 0; blue, 0 }  ][line width=0.75]    (10.93,-3.29) .. controls (6.95,-1.4) and (3.31,-0.3) .. (0,0) .. controls (3.31,0.3) and (6.95,1.4) .. (10.93,3.29)   ;

	\draw    (333,36) -- (350,36) ;
	\draw [shift={(352.54,36)}, rotate = 179.4] [color={rgb, 255:red, 0; green, 0; blue, 0 }  ][line width=0.75]    (10.93,-3.29) .. controls (6.95,-1.4) and (3.31,-0.3) .. (0,0) .. controls (3.31,0.3) and (6.95,1.4) .. (10.93,3.29)   ;

	\draw  [fill={rgb, 255:red, 155; green, 155; blue, 155 }  ,fill opacity=1 ] (487.36,162.22) -- (437.55,209.96) -- (437.55,98.22) -- (487.36,50.48) -- cycle ;

	\draw  [fill={rgb, 255:red, 155; green, 155; blue, 155 }  ,fill opacity=1 ] (437.55,99.18) -- (484.9,51.83) -- (598.92,51.83) -- (598.92,162.32) -- (551.57,209.67) -- (437.55,209.67) -- cycle ; \draw   (598.92,51.83) -- (551.57,99.18) -- (437.55,99.18) ; \draw   (551.57,99.18) -- (551.57,209.67) ;

	\draw    (471.73,107.6) .. controls (510.66,99.63) and (511.35,39.15) .. (514.77,27.69) ;

	\draw    (475.15,139.53) .. controls (527.75,132.39) and (544.68,219.3) .. (542.12,224.22) ;

	\draw    (471.73,113.33) .. controls (513.22,107.82) and (537.85,68.63) .. (543.83,30.14) ;

	\draw    (476.86,147.72) .. controls (517.02,146.9) and (517.34,235.68) .. (514.77,225.85) ;

	\draw    (477.71,130.53) .. controls (545.22,120.7) and (574.59,220.94) .. (572.03,225.85) ;

	\draw    (472.97,121.3) .. controls (514.46,115.79) and (569.46,69.45) .. (575.45,30.96) ;

	\draw  [dash pattern={on 4.5pt off 4.5pt}] (487.36,50.49) -- (487.36,165.98) -- (437.55,209.67) ;

	\draw (140,27) node [anchor=north west][inner sep=0.75pt]   [align=left] {x};

	\draw (13,235) node [anchor=north west][inner sep=0.75pt]   [align=left] {y};

	\draw (43,106) node [anchor=north west][inner sep=0.75pt]   [align=left] {$\displaystyle \Gamma _{0}$};

	\draw (410,27) node [anchor=north west][inner sep=0.75pt]   [align=left] {x};

	\draw (11,6.73) node [anchor=north west][inner sep=0.75pt]   [align=left] {Case 1};

	\draw (292,2.37) node [anchor=north west][inner sep=0.75pt]   [align=left] {Case 2};

	\draw (275,235) node [anchor=north west][inner sep=0.75pt]   [align=left] {y};

	\draw (76,25) node [anchor=north west][inner sep=0.75pt]   [align=left] {$\displaystyle \phi |_{\Omega_0}$};

	\draw (15,25) node [anchor=north west][inner sep=0.75pt]   [align=left] {$\displaystyle \phi |_{\Gamma _0}$};

	\draw (134,105) node [anchor=north west][inner sep=0.75pt]   [align=left] {\small TIE};

	\draw (291,25) node [anchor=north west][inner sep=0.75pt]   [align=left] {$\displaystyle \phi |_{\Omega_0}$};

	\draw (86,122) node [anchor=north west][inner sep=0.75pt]   [align=left] {$\displaystyle \Omega _{0}$};

	\draw (195,121) node [anchor=north west][inner sep=0.75pt]   [align=left] {$\displaystyle \Omega _{0}$};

	\draw (353,123) node [anchor=north west][inner sep=0.75pt]   [align=left] {$\displaystyle \Omega _{0}$};

	\draw (358,25) node [anchor=north west][inner sep=0.75pt]   [align=left] {$\displaystyle \phi |_{\Omega }$};

	\draw (399,105) node [anchor=north west][inner sep=0.75pt]   [align=left] {\small TPE};

	\draw (607.36,141.32) node [anchor=north west][inner sep=0.75pt]   [align=left] {z};

	\draw (449.08,135.22) node [anchor=north west][inner sep=0.75pt]   [align=left] {$\displaystyle \Omega _{0}$};

\end{tikzpicture}
    \caption{Two boundary cases for the TIE and the TPE.}
    \label{fig:introboundary}
\end{figure}

\textbf{Outline.} In the upcoming section, we derive the TIE and the TPE. \autoref{sec:TIE} is dedicated to solving the TIE numerically for two different phase functions $\phi$. Following that, in \autoref{sec:TPE}, we discuss the TPE. We investigate local existence of the TPE using the method of characteristics (see \autoref{subsection: charac}) and we study viscosity solutions of the TPE in \autoref{subsection: visco}. All presented numerical simulations were implemented with MATLAB. Finally, in \autoref{sec:comparison} we consolidate our findings and discuss the results from \autoref{sec:TIE} and \autoref{sec:TPE}.

\subsection*{Notation}
\label{sec:notation}
Consider a function $f:\R^3 \times [0,\infty)$, $(x,y,z,t) \mapsto f(x,y,z,t)$.
\begin{itemize}
	\item We denote its partial derivatives by $f_x,f_y,f_z,f_t$,
	\item $\nabla f= (f_x,f_y,f_z)^T$ is the spatial gradient of $f$.
	\item The $n$-th partial derivatives are denoted by $f_{\rho_1,\cdots,\rho_n}$, where $\rho_i \in \set{x,y,z,t}$.
	\item $\Delta f= f_{xx} + f_{yy} + f_{zz}$ denotes the three-dimensional Laplace operator in space. 
	\item $\vx=(x,y)$ denotes a vector in the $xy$-plane of the $xyz$-space $\R^3$.
	\item $\nabla_{\vx} f =(f_x,f_y)^T$ denotes the gradient of $f$ with respect to $\vx$.
	\item The Laplace operator with respect to $\vx$ is denoted by $\Delta_{\vx}f = f_{xx} + f_{yy}$.
\end{itemize}

\section{Derivation of the Transport of Intensity and Phase Equations}
\label{sec:derivation}
In this section we carry out the derivation of the transport of intensity (TIE) and the transport of phase equation (TPE). For this purpose we derive the paraxial Helmholtz equation first, where we closely follow \cite{GurRobNug95}, the details are taken from \cite{Goodman} and \cite{Saleh2019}.

\subsection{Paraxial Helmholtz equation} 
We derive the mathematical equation for a wave $u$ traveling preferably in $z$ direction, the so-called paraxial Helmholtz equation. The starting point is the \emph{scalar wave equation}
\begin{equation}\label{eq:wave}
	\frac{n^2}{c^2} u_{tt}(\vx,z,t) - \Delta u(\vx,z,t) =0,	
\end{equation}
where $n$ represents the refractive index of the medium and $c$ denotes the vacuum velocity of light.

A monochromatic wave with angular velocity $\omega_0$ is given by
\begin{equation} \label{eq:timeharmonic}
	u(\vx,z,t)= u_0(\vx,z) \e^{\i\omega_0t}.
\end{equation}
If  $u_0$ solves the Helmholtz equation 
 \begin{equation}\label{eq:helmholtz}
 	\Delta u_0 (\vx,z)+k^2 u_0(\vx,z) =0
 \end{equation}
with wave number $k:= \frac{n \omega_0}{c}$, then the monochromatic wave $u$ solves \autoref{eq:wave}.

We follow \cite{Saleh2019} to derive the \emph{paraxial Helmholtz equation} by making the ansatz 
\begin{equation}\label{eq:paraxansatz}
	u_0(\vx,z)=A(\vx,z)\e^{-\i kz}
\end{equation}
and construct $A$ in such a way that $u_0$ solves \autoref{eq:helmholtz}, which means that
\begin{equation*}\label{eq:au0} \begin{aligned}
	0 &= \Delta_{\vx} \left(A(\vx,z)\e^{-\i kz}\right) +\frac{\partial ^2}{\partial z^2}\left(A(\vx,z)\e^{-\i kz}\right) +k^2 A(\vx,z)\e^{-\i kz}\\
	&=\e^{-\i kz}\left(\Delta_{\vx} A + A_{zz} -2 \i kA_{z}-k^2A(\vx,z)\right) + k^2A(\vx,z)\e^{-\i kz}.
\end{aligned} \end{equation*}
After division by the factor $\e^{-ikz}$, we see that $A$ is a solution of  
\begin{equation} \label{eq:Afull}
\Delta_{\vx} A + A_{zz} -2 \i kA_{z}= 0. 
\end{equation}
Next, we assume that relative variations of $A$ at scale distance $\lambda$ are small in comparison with the wavelength $k = \frac{2\pi}{\lambda}$, which means that
\begin{equation}\label{eq:parax1}
	\frac{ \abs{ A_{z}(\vx,z)}} {\abs{A(\vx,z)}} \sim \frac{1}{\lambda} 
    \frac{ \abs{A(\vx,z+\lambda)-A(\vx,z)}}{\abs{A(\vx,z)}} \ll k .
\end{equation}
Note that all dependencies of $A$ are complex. We make the same assumption for $A_{z}$, that is 
\begin{equation}\label{eq:parax2}
	\frac{\abs{A_{zz}(\vx,z)}}{\abs{A_{z}(\vx,z)}} \ll k,
\end{equation}
and then from \autoref{eq:parax2} and \autoref{eq:parax1} we deduce that
\begin{equation} \label{eq:A}
	\Delta_{\vx} A -2\i kA_{z}= 0,
\end{equation}
which is the \emph{paraxial Helmholtz equation}.

\subsection{Derivation of the TIE and the TPE from the paraxial Helmholtz equation}
We use the representation of a complex wave $A$ as in \autoref{eq:AI}. 

With this representation it follows that, 
\begin{equation*}	\label{eq:dervergleich}
	2 \i kA_{z} 
	=k \e^{\i\phi }\left(-2\sqrt{I}\phi_z+\i \frac{I_{z}}{\sqrt{I}}\right)
\end{equation*}
and
\begin{equation} \label{eq:tietpe}
	\Delta_{\vx}A =\e^{\i\phi }\left( -\frac{\| \nabla_{\vx}I\|^2}{4I^{\frac{3}{2}}}+\frac{\Delta_{\vx}I}{2\sqrt{I}}
	-\sqrt{I}\| \nabla_{\vx}\phi\|^2+\i\sqrt{I}\Delta_{\vx}\phi 
	+\i\frac{\nabla_\vx I \cdot \nabla_\vx \phi}{\sqrt{I}}\right),
\end{equation}
where we left out the spatial variables $(\vx,z)$ for notational convenience.
Inserting these expressions into the 
paraxial Helmholtz \autoref{eq:A} and dividing by $\e^{\i\phi}$, we get two equations, one for the real and one for the imaginary part. For the equation resulting from the imaginary part we obtain by multiplication with $\sqrt{I}$ the TIE (see \autoref{eq:TIE}).
Dividing the real part of equation \autoref{eq:tietpe} by $\sqrt{I}$ we get the TPE (see \autoref{eq:TPE}).

\subsection{Modelling errors of the TIE and the TPE}
Both the TIE and the TPE relate the intensity $I$ to the phase $\phi$. Therefore, the question arises, whether one of these equations is more accurate, in the sense that the left out second order term in the paraxial approximation affects one of the terms (imaginary part or real part) less than the other. When calculating the second derivative $A_{zz}$, we obtain
\begin{equation*}
	A_{zz} = \e^{\i \phi} 
	\underbrace{\left(  -\frac{1}{4} \frac{I_{z}^2}{I^\frac{3}{2}}  + \frac{1}{2} \frac{I_{zz}}{\sqrt{I}}  - 
		 \sqrt{I} \phi_z^2 \right)}_{=:M_{TPE}}
	+ \i \e^{\i \phi}\underbrace{\left(\frac{I_{z}}{\sqrt{I}} \phi_z + \sqrt{I} \phi_{zz}\right)}_{=:M_{TIE}}.
\end{equation*}
The real and imaginary parts that are omitted describe the modeling errors of the TPE and the TIE, respectively.
We see that for the TIE first and second order derivatives of the phase are lost in contrary to the TPE, where only the first order term is lost. However, if one assumes for instance constant intensity of $I$ in the $z$-direction then the first order term in the TIE loss vanishes, whereas the corresponding lost term for the TPE is still present. 
In summary, there is no immediate sign that either one of the two equations TIE or TPE is more accurate. Since the TIE is linear for given $I$ and the TPE is nonlinear there is a preference for the TIE, which is why it is used in the literature often \cite{BosFroNilSagUns16, Petrucelli}.

\section{The Transport of Intensity Equation (TIE)}
\label{sec:TIE}
The TIE and the TPE, in \autoref{eq:TIE} and \autoref{eq:TPE}, respectively, were derived in free space $\R^3$. For practical computations, however, boundary information on the phase is mandatory and without educated guesses or measurements on the phase at the boundary, the recovered phase information inside $\Omega$ cannot be accurate. 

In what follows we always assume that the basic assumptions of the theory of elliptic partial differential equations are satisfied (see \cite{Hac18}). In particular, this means that the intensity $I$ of the wave is always positive and the domain is regular. In \cite{Adams} different types of regularity are discussed, for a detailed discussion of the matter in the context of proving existence and uniqueness results of weak solutions of the TIE we refer to \cite{Str23}. 

In the following, we will solve the TIE on $\Omega_0 \subseteq \R^2$ given some Dirichlet boundary condition $g$ on $\Gamma_0$. This corresponds to the first case of the possible settings introduced in \autoref{sec:introduction}.
We will consider the problem 
\begin{equation*}\label{eq:stronginhomdir}	\nabla_\vx\cdot (I\nabla_\vx\phi)=k I_{z} \quad \text{ in } \Omega_0 \quad \text{ and } \quad 
    \phi \vert_{\Gamma_0}=g, 
\end{equation*}
for two different ground-truth functions $A = \sqrt{I}\e^{\i\phi}$ solving the paraxial Helmholtz equation and we analyze the performance of the TIE for reconstructing the phase $\phi$ given $I$.

\subsection{Numerical example: constant intensity $(I \equiv 1)$}
\label{subsec:ex1}
\begin{figure}[h]
    \centering
    \begin{minipage}{.32\textwidth}
    \begin{tikzpicture}
            \node (img) {\includegraphics[width=\textwidth]{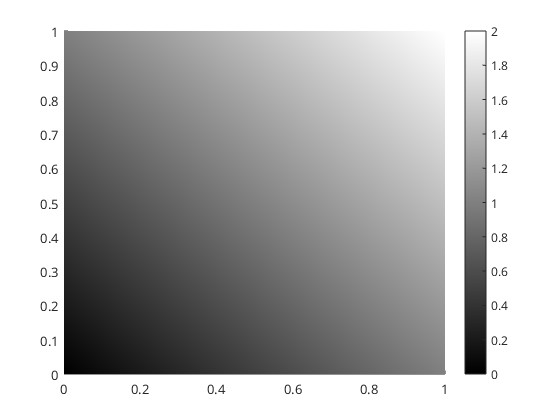}};
            \node[below=of img,yshift=1.3cm,xshift=-0.2cm]{\tiny$x$};
                                    \node[left=of img, rotate=90,anchor=center,yshift=-1.2cm]{\tiny $y$};
    \end{tikzpicture}
            \subcaption{Ground-truth}
            \label{fig:pw_groundtruth}
    \end{minipage}
    \begin{minipage}{.32\textwidth}
    \begin{tikzpicture}
            \node (img){ \includegraphics[width=\textwidth]{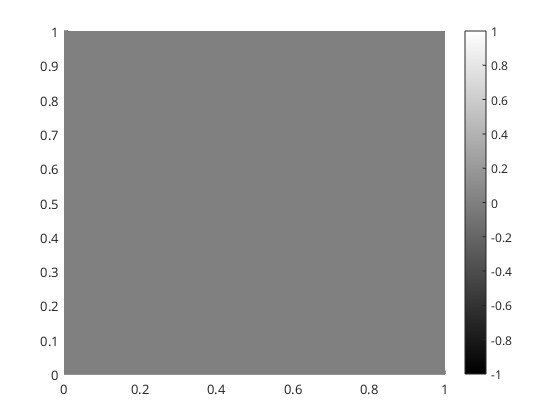}};
                        \node[below=of img,yshift=1.3cm,xshift=-0.2cm]{\tiny$x$};
                        \node[left=of img, rotate=90,anchor=center,yshift=-1.2cm]{\tiny $y$};
                \end{tikzpicture}
            \subcaption{$\phi\vert_{\Gamma_0}=0$}
            \label{fig:pw_zerodirichlet}
    \end{minipage}
    \begin{minipage}{.32\textwidth}
    \begin{tikzpicture}
        \node (img) { \includegraphics[width=\textwidth]{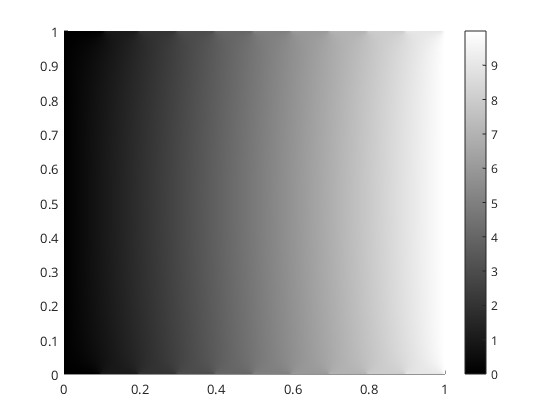}};
                    \node[below=of img,yshift=1.3cm,xshift=-0.2cm]{\tiny$x$};
                                            \node[left=of img, rotate=90,anchor=center,yshift=-1.2cm]{\tiny $y$};
    \end{tikzpicture}
         \subcaption{$\phi\vert_{\Gamma_0}= \lfloor 10 x \rfloor$}
         \label{fig:pw_floorbc}
    \end{minipage}
    \caption{Ground-truth of the phase $\phi(\vx, z=0)$ described in \autoref{eq:phase:ex}, (left), reconstruction with $\phi\vert_{\Gamma_0}=0$ (middle) and $\phi\vert_{\Gamma_0}= \lfloor 10 x \rfloor$ (right)}
    \label{fig:phase:ex}
\end{figure}

\begin{figure}[h]
    \centering
    \begin{minipage}{0.4\textwidth}
    \begin{tikzpicture}
        \node (img) {        \includegraphics[width=\textwidth]{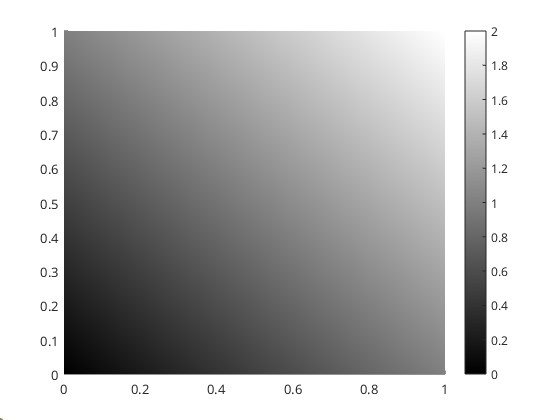}};
             \node[below=of img,yshift=1.3cm,xshift=-0.2cm]{\tiny$x$};
                                            \node[left=of img, rotate=90,anchor=center,yshift=-1.2cm]{\tiny $y$};
    \end{tikzpicture}
        \subcaption{Reconstruction}
        \label{fig:pw_reconstruction_correct}
    \end{minipage}
     \begin{minipage}{0.41\textwidth}
     \begin{tikzpicture}
         \node (img) {        \includegraphics[width=\textwidth]{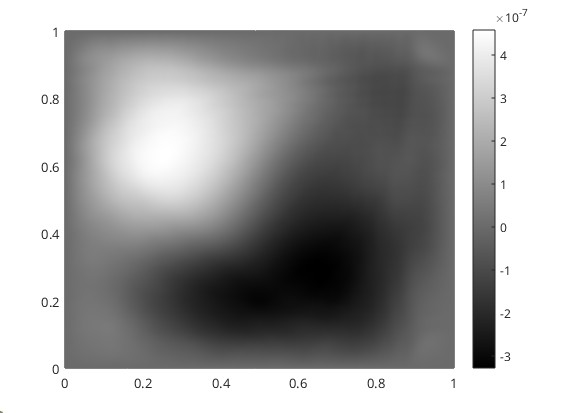}};
              \node[below=of img,yshift=1.3cm,xshift=-0.2cm]{\tiny$x$};
                                            \node[left=of img, rotate=90,anchor=center,yshift=-1.2cm]{\tiny $y$};
     \end{tikzpicture}
        \subcaption{Reconstruction error}
        \label{fig:pw_reconstruction_error}
    \end{minipage}
    \caption{Reconstruction (left) and reconstruction error (right) of the phase $\phi(\vx, z=0)$ described in \autoref{eq:phase:ex} given the ground-truth boundary condition}
    \label{fig:TIE:plane_wave}
\end{figure}
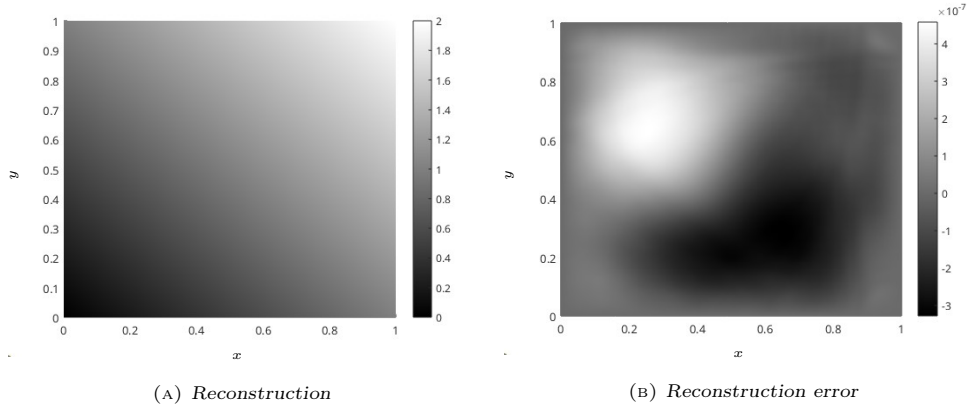

Let $\vek \xi \in \R^2$ be fixed. Then, a short calculation shows that the function
\begin{align*} 
	A(\vx, z)=  \e^{\i \left(\vx \cdot \vek \xi+  \frac{1}{2} \| \vek \xi \|^2 z\right)} \text{ for } \vx \in \R^2, z \in \R
\end{align*}
is a solution of the paraxial Helmholtz \autoref{eq:A} with wavenumber $k=1$ (see \cite{Eva10b}). Comparing with the general form of $A$ in \autoref{eq:AI} we see that 
\begin{align} \label{eq:phase:ex}
	I(\vx, z) = 1 \text{ and } 
	\phi(\vx,z) = \vx \cdot \vek \xi+  \frac{1}{2} \| \vek \xi \|^2 z.
\end{align}
\autoref{fig:pw_groundtruth} depicts the function $\vx \in \R^2 \to \phi(\vx,0)$ with $\vek \xi=(1,1)^T$ as specified in \autoref{eq:phase:ex}.
Since in this case the intensity is constant, the TIE \autoref{eq:TIE} reduces to the Laplace equation for $\phi$, which we solved for different boundary conditions. The numerical solution was obtained by means of a finite element method with linear interpolation functions. A mesh on the square $[0,1] \times [0,1]$ with 11821 nodes was created using Gmsh \cite{gmsh}.  We observe that the boundary data dominates the outcome and in general we will not be able to retrieve the phase, if we do not impose the ground-truth boundary data given by \autoref{eq:phase:ex} as input data for the reconstruction. The corresponding results for $\phi\vert_{\Gamma_0}=0$ and $ \phi\vert_{\Gamma_0}=\lfloor 10x \rfloor$ are depicted in \autoref{fig:pw_zerodirichlet} and \autoref{fig:pw_floorbc}, respectively. Moreover, the reconstruction of the phase in $\Omega_0$ with the ground-truth boundary values is shown in \autoref{fig:pw_reconstruction_correct}. In addition, we also plot the error that is made in the reconstruction (see \autoref{fig:pw_reconstruction_error}).

\subsection{Numerical example: constant phase $\phi \equiv \frac{3\pi}{2}$}
\label{subsec:ex2}
We have seen that in the case that the intensity is constant but the phase is not, we can reconstruct the phase reasonably well if we impose the boundary condition of the ground-truth on the TIE. Now we consider a numerical example where the phase is constant in the plane $\{z=0\}$, whereas the intensity is not. The Gaussian beam (see \cite{Saleh2019}) is defined by
\begin{align}
   A(\vx,z)= I_0\frac{1}{q(z)}\e^{-\i k \frac{\norm{\vx}^2}{2q(z)}},
   \label{eq:gbeam}
\end{align}
where
\begin{equation*}
    q(z):= z + \i z_R
\end{equation*}
and $z_R$ denotes the so-called Rayleigh range. For our example we consider $I_0=1$.
It serves as an ideal example for testing the reconstruction of the phase by both the TIE and the TPE, since it fulfills the paraxial Helmholtz equation (see \cite{Saleh2019}).
We define further
\begin{equation*}
        w(z):= w_0 \sqrt{1+ \left( \frac{z}{z_R}\right)^2 },
\end{equation*}
where 
\[
w_0=\sqrt{\frac{\lambda z_R}{\pi}}
\]

 is the waist radius of the Gaussian beam, see \cite{Saleh2019} for a detailed description.

Furthermore, it holds that
\begin{equation*}
    \frac{1}{q(z)} = \frac{1}{R(z)}-\i \frac{\lambda}{\pi w^2(z)}
\text{ with }
        R(z):= z\left(1+ \left( \frac{z_R}{z}\right)^2\right),
\end{equation*}
such that we can rewrite \autoref{eq:gbeam} as
\begin{align}
    A(\vx,z) = \frac{w_0}{z_Rw(z)} \e^{-\frac{\norm{\vx}^2}{w(z)^2}}\e^{-\frac{\i k\norm{\vx}^2}{2R(z)}} \e^{\i \operatorname{arg}\left( \frac{1}{q(z)}\right)},
    \label{eq:gaussianbeam_seperated}
\end{align}
where $\operatorname{arg}\left(\frac{1}{q(z)}\right)$ denotes the argument of the complex number $\frac{1}{q(z)}$.
In particular, we obtain for $z=0$ that
\begin{align}\label{eq:i_intvarphicost}
    I(\vx,z=0) = \frac{1}{z_R^2}\e^{-2\frac{\norm{\vx}^2}{w_0^2}} \text{ and } 
    \phi(\vx,z=0) = \operatorname{arg} \left( \frac{1}{\i z_R} \right),
\end{align}
which shows that the phase is constant, while the intensity is not.
In our experimental setup we choose $k=1$ and $z_R=1$, such that $w_0= \sqrt{\frac{2 z_R}{k}}= \sqrt{2}$, which implies that
\begin{equation}\label{eq:gbphi}
    I(\vx, z=0)= \e^{-\norm{\vx}^2} \text{ and } 
    \phi(\vx,z=0) = \operatorname{arg} \left( \frac{1}{\i} \right) = \frac{3 \pi}{2}.
\end{equation}
For homogeneous Dirichlet-boundary conditions we again obtain the zero solution.

\begin{figure}[h]
    \centering
    \begin{minipage}{0.4\textwidth}
    \begin{tikzpicture}
        \node (img) {		\includegraphics[width=\textwidth]{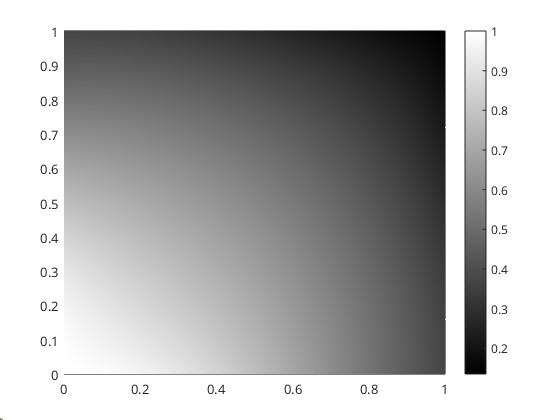}};
                  \node[below=of img,yshift=1.3cm,xshift=-0.2cm]{\tiny$x$};
                                            \node[left=of img, rotate=90,anchor=center,yshift=-1.2cm]{\tiny $y$};
    \end{tikzpicture}
        \subcaption{$\phi\vert_{\Gamma_0}= \e^{-\norm{\vx}^2}$}
        \label{fig:gb_random1}
    \end{minipage}
     \begin{minipage}{0.4\textwidth}
     \begin{tikzpicture}
         \node (img) {		\includegraphics[width=\textwidth]{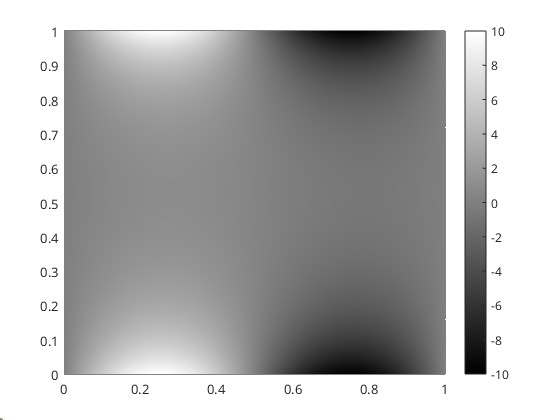}};
                   \node[below=of img,yshift=1.3cm,xshift=-0.2cm]{\tiny$x$};
                                            \node[left=of img, rotate=90,anchor=center,yshift=-1.2cm]{\tiny $y$};
     \end{tikzpicture}
        \subcaption{$\phi\vert_{\Gamma_0}=10\sin(2\pi x)$}
        \label{fig:gb_random2}
    \end{minipage}
\caption{Reconstruction of the phase $\phi(\vx, z=0)$ described in \autoref{eq:gbphi} with boundary data not corresponding to the ground-truth boundary value.
	}     \label{fig:TIE_numerical2}
\end{figure}
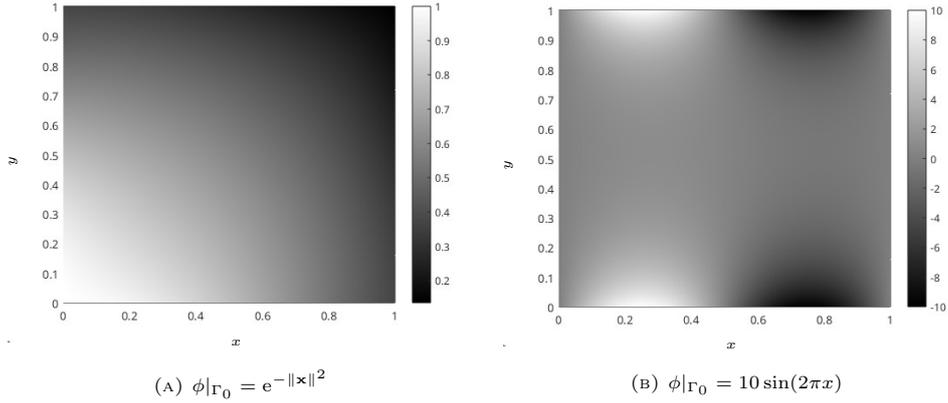
\begin{figure}[h]
    \centering
    \begin{minipage}{0.4\textwidth}
    \begin{tikzpicture}
        \node (img) {            \includegraphics[width=\textwidth]{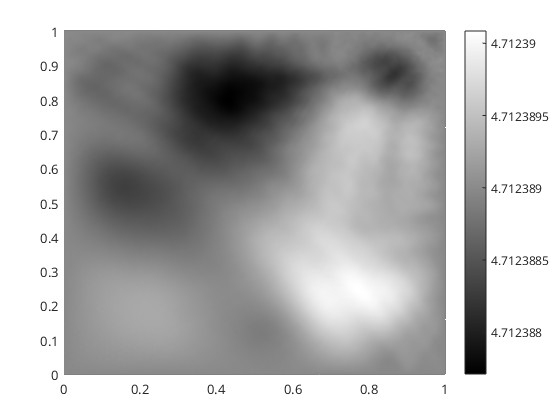}};
                  \node[below=of img,yshift=1.3cm,xshift=-0.2cm]{\tiny$x$};
                                            \node[left=of img, rotate=90,anchor=center,yshift=-1.2cm]{\tiny $y$};
    \end{tikzpicture}
            \subcaption{Reconstruction}
            \label{fig:gb_correct_reconstruction}
    \end{minipage}
       \begin{minipage}{0.41\textwidth}
       \begin{tikzpicture}
           \node (img) {    \includegraphics[width=\textwidth]{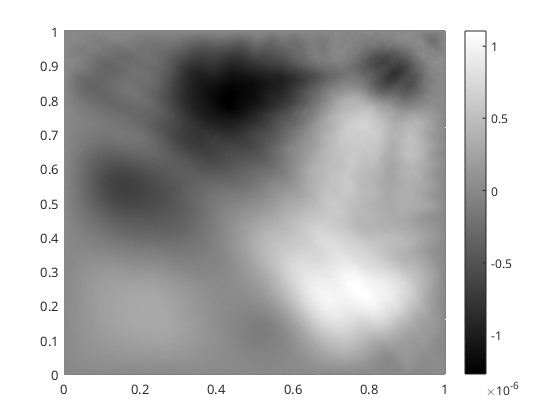}};
                     \node[below=of img,yshift=1.3cm,xshift=-0.2cm]{\tiny$x$};
                                            \node[left=of img, rotate=90,anchor=center,yshift=-1.2cm]{\tiny $y$};
       \end{tikzpicture}
            \subcaption{Reconstruction Error}
            \label{fig:gb_error}
    \end{minipage}
    \caption{Reconstruction (left) and reconstruction error (right)  of the phase $\phi(\vx, z=0)$ described in \autoref{eq:gbphi} given the ground-truth boundary condition}
    \label{fig:TIE:gaussian}
\end{figure}
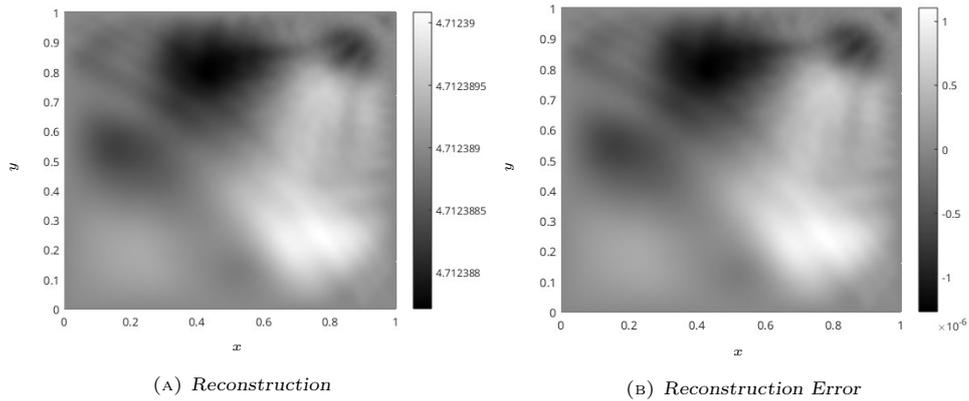

\autoref{fig:TIE_numerical2} shows the solutions in the case that we impose boundary conditions that do not correspond to the known ground-truth boundary data. The results for the boundary value $\e^{-\norm{\vx}^2}$ and $10\sin(2\pi x)$ are depcited in \autoref{fig:gb_random1} and \autoref{fig:gb_random2} respectively. 
Again, we compute the phase given the ground-truth boundary data: The result is depicted in \autoref{fig:gb_correct_reconstruction} and we see that the phase is reconstructed reasonably well.
Additionally, also plotting the error (see figure \autoref{fig:gb_error}) confirms this.

\section{Solving the Transport of Phase Equation}
\label{sec:TPE}
In this section, we investigate methods for solving the nonlinear TPE \autoref{eq:TPE} with the method of characteristics and the viscosity method. With the first one we can obtain analytical solutions for special cases. The viscosity solutions provide approximate solutions, which can be calculated with finite element methods. 
 
\subsection{Characteristic method} \label{subsection: charac}
The TIE allows us to compute the phase $\phi$ in $\Omega_0$ from phase data on $\Gamma_0$ (see \autoref{subsec:ex1}). 
Now we investigate the method of characteristics to determine local solutions of the TPE. 
See \cite[Chap.\ 3.2]{Eva10b}.
We aim at finding a solution of \autoref{eq:TPE} on curves starting from points $(\vx,0) \in \Omega_0$ given data
\begin{equation}
\phi(\vx, 0) = g(\vx), \quad 0 \le x,y \le 1.
	\label{eq:characTPE} 
\end{equation}
The TPE is a nonlinear partial differential equation of first order for a function $\phi:\R^3 \to \R$. We first write it as
\begin{equation}\label{eq:defFcharact}
	F(\phi_x, \phi_y, \phi_z, \phi, x, y, z) = 0, 
\end{equation}
where the function $F$ is given by
\begin{equation*}
		F:\R^7 \to \R, \quad  (p,q,r,s,x,y,z) \mapsto 2kr -(p^2+q^2) +\hat{I}(x,y,z).
\end{equation*}
Then we introduce the characteristics, that is, the arguments of $F$ parametrized with respect to $\tau \in \R$, that is, 
\begin{equation*}
	\begin{aligned}
		\tau &\mapsto (x(\tau),y(\tau),z(\tau)), \\
		\tau &\mapsto s(\tau) = \phi (x(\tau),y(\tau),z(\tau)),\\
		\tau &\mapsto p(\tau) = \phi_x (x(\tau),y(\tau),z(\tau)),\\ 
		\tau &\mapsto q(\tau) = \phi_y (x(\tau),y(\tau),z(\tau)),\\
		\tau &\mapsto r(\tau) = \phi_z (x(\tau),y(\tau),z(\tau)).
	\end{aligned}
\end{equation*}
After parametrization, the partial differential equation can be solved by computing the solution of the following \emph{characteristic ODEs}. Below, differentiation with respect to $\tau$ is abbreviated by a prime $(\cdot)'$.
\begin{equation}
	\begin{split}
		x'	&= \partial_p F = -2p, \qquad y'  = \partial_q F = -2q, \qquad z'  = \partial_r F = 2k, \\
		s'	&= p\partial_p F + q\partial_q F + r\partial_r F = -2p^2 -2q^2 + 2kr, \\
		p'	&= -\partial_x F -p\partial_sF = -\hat{I}_x\\
			&= -\frac{(I_x^2+I_y^2)I_x}{2I^3}+\frac{I_{xx}I_x+I_{yx}I_y}{2I^2}+\frac{I_x(I_{xx}+I_{yy})}{2I^2}-\frac{I_{xxx}+I_{yyx}}{2I},\\
		q'  &= -\partial_y F -q\partial_sF = -\hat{I}_y\\
			&= -\frac{(I_x^2+I_y^2)I_y}{2I^3}+\frac{I_{xy}I_x+I_{yy}I_y}{2I^2}+\frac{I_y(I_{xx}+I_{yy})}{2I^2}-\frac{I_{xxy}+I_{yyy}}{2I},\\
		r'	&= -\partial_z F -r\partial_sF = -\hat{I}_z\\
			&= -\frac{I_z(I_x^2+I_y^2)}{2I^3}+\frac{I_xI_{xz}+I_yI_{yz}}{2I^2}+\frac{I_z(I_{xx}+I_{yy})}{2I^2}-\frac{I_{xxz}+I_{yyz}}{2I}.\\
	\end{split}
	\label{eq:characteristic equations}
\end{equation}
For an arbitrary starting point $(\vx_0,0) = (x_0,y_0,0) \in \Omega_0$ we complement the characteristic equations with the following initial conditions, also called compatibility conditions,
\begin{equation*}
	\begin{split}
		x(0) &= x_0, \quad y(0) = y_0, \quad z(0) = 0,\\
		s(0) & = g(x_0,y_0) \eqqcolon s_0,\\
		p(0) &= g_x(x_0,y_0) \eqqcolon p_0,\\
		q(0) &= g_y(x_0,y_0)\eqqcolon q_0,\\
		r(0) &= \frac{p_0^2 + q_0^2 - \hat I (x_0,y_0,0)}{2k} \eqqcolon r_0.
	\end{split}
\end{equation*}
Note that $\partial F/ \partial r = 2k$ and that the wave number $k$ is positive. Thus, the noncharacteristic condition
\begin{equation*}
	\frac{\partial F}{\partial r}(p_0, q_0, r_0, s_0, x_0, y_0, 0) \neq 0
\end{equation*}
is satisfied, guaranteeing existence of a local solution of the TPE close to $(\vx_0,0)$.

%

\subsubsection{Example: constant intensity in $xy$-plane}
Assume that $I_x = I_y = 0$ in $\Omega$. Then $\hat I = 0$  and \autoref{eq:characteristic equations} implies that
\begin{equation*}
	\begin{split}
		p'	= q'   =  r'   = 0.
	\end{split}	
\end{equation*}
Thus the solutions of the characteristic equations are given by 
\begin{equation}
	\begin{split}
		x(\tau) &= -2 p_0 \tau + x_0, \quad y(\tau) = -2 q_0 \tau + y_0, \quad z(\tau) = 2k\tau,\\
		s(\tau)	&= \left(-2p_0^2 - 2q_0^2 + 2kr_0 \right)\tau + s_0 = -\left(p_0^2 + q_0^2 \right)\tau + s_0, \\
		p(\tau) &= p_0, \quad q(\tau) = q_0, \quad r(\tau) = r_0.
	\end{split}
	\label{eq:charauniforsolution}
\end{equation}
For every point $(\vx_0 ,0)$, the characteristic $\tau \mapsto (x(\tau),y(\tau),z(\tau))$ passing through it is a line. If $g$ is of the form $g(x,y) = \alpha x + \beta y + \gamma$, for $\alpha, \beta, \gamma \in \R$, then these lines are all parallel and we obtain a global solution of the TPE
\begin{equation*}
	\phi(x,y,z) = g(x,y) + z\frac{\alpha^2 + \beta^2}{2k}, \quad (x,y,z) \in \R^3, 
\end{equation*}
which satisfies \autoref{eq:characTPE} on the entire plane $\{z=0\}$.


\subsubsection{Example: Gaussian beam}
The intensity of the Gaussian beam, as defined in \autoref{subsec:ex2}, is given by  
\begin{equation*} 
I(\vx, z) = \frac{w_0^2}{w^2(z)} \e^{-2 \frac{\norm{\vx}^2}{w^2(z)}}.
\end{equation*} 
The ODE system \autoref{eq:characteristic equations} then becomes
\begin{equation*}
	\begin{split}
		x'  &= -2p, \quad y'  = -2q, \quad z'  = 2k,\\
		s'	&= -2p^2 -2q^2 + 2kr,\\
		p'  &= -8xw^{-4}(z), \quad q'  = -8yw^{-4}(z),\\
		r' 	&= -\frac{8 z z_R^2 (-2 (x^2 + y^2) z_R^2 + w_0^2 (z^2 + z_R^2))}{w_0^4 (z^2 + z_R^2)^3}.
	\end{split}
	\label{eq:gaussianbeamcharac}
\end{equation*}
However, it is difficult to solve this ODE system. The characteristic method cannot provide an explicit solution for the Gaussian beam. In the following we use the viscosity method to find an approximate solution of the TPE. 

\subsection{Viscosity solutions of the TPE} \label{subsection: visco}
We aim to solve the TPE, \autoref{eq:TPE}, in $\Omega = [0,1]^3$ together with an initial condition $\phi(\vx,0)=g(\vx)$. Our numerical solution is based on a viscosity approximation, which requires a virtual boundary condition on $\Gamma = \bigcup_{0\le z \le 1} \Gamma_z$ and a parameter $\epsilon >0$.
The viscosity approximation of the TPE on $\Omega$ is then given by:
\begin{equation}\label{eq:visc-bound}
	\begin{dcases}
		2k \phi_z^\epsilon - \epsilon \Delta_\vx\phi^\epsilon - \|\nabla_{\vx} \phi^\epsilon\|^2 = - \hat{I}, & \text{in } \Omega,	\\
			\phi^\epsilon = g, & \text{on } \Omega_0, \\
			\phi_\epsilon = h, & \text{on } \Gamma.
	\end{dcases}
\end{equation}
If the limit of $\phi^\epsilon$ for $\epsilon \to 0+$ exists, then it is considered a generalized solution of the TPE and it is named the \textit{vanishing viscosity solution} \cite{CraLio83b}.
Therefore, the fully nonlinear TPE is approximated by a quasilinear parabolic partial differential equation \autoref{eq:visc-bound}. The additional term $\epsilon\Delta\phi^\epsilon$ acts as a regularizer.

\subsubsection{Solution of the quasilinear PDE}
We use the Cole-Hopf transformation (see \cite[p.\ 206]{Eva10b}) to transform \autoref{eq:visc-bound} to a linear equation, which can be solved with standard numerical methods. The Cole-Hopf transformation $\psi^{\epsilon}$ of $\phi^\epsilon$ is defined by the relation
	\[\phi^\epsilon = \epsilon  \log \psi^{\epsilon}.\]
It satisfies 
\begin{equation}\label{eq:colehopf}
	\begin{dcases}
		2k\psi_z^{\epsilon} - \epsilon \Delta_\vx\psi^{\epsilon} = -\frac{\hat{I}}{\epsilon}\psi^{\epsilon}, & \text{ in } \Omega, \\
		\psi^{\epsilon} = \exp\left(\frac{g}{\epsilon}\right), & \text{ on } \Omega_0, \\
        \psi^{\epsilon} = \exp\left(\frac{h}{\epsilon}\right), & \text{ on } \Gamma,
	\end{dcases}
\end{equation}
where we have used that
\[ \phi^{\epsilon}_z = \epsilon \frac{\psi^{\epsilon}_z}{\psi^{\epsilon}}, \qquad
	\Delta_\vx \phi^\epsilon = \frac{\epsilon \Delta_\vx \psi^{\epsilon}}{\psi^{\epsilon}}- \frac{\epsilon \norm{\nabla_\vx \psi^{\epsilon}}^2}{(\psi^{\epsilon})^2}, \qquad
	\nabla_\vx \phi^\epsilon = \frac{\epsilon \nabla_\vx \psi^{\epsilon}}{\psi^{\epsilon}}.\]
We consider an example similar to the one for the TIE equation: 
\begin{example}\label{ex:viscosity_gauss}
We assume a Gaussian beam intensity
\[I(\vx, z) = \frac{w^2_0}{w^2(z)}\exp\left(-2\frac{\norm{\vx}^2}{w\left(z\right)^2}\right) \text{ for } (\vx,z) \in \Omega\]
as in \autoref{eq:gaussianbeam_seperated}. The approximate phase $\phi^\epsilon$ is computed by first solving \autoref{eq:colehopf} for $\psi^\epsilon$ and applying the Cole-Hopf transformation later. To be able to compare with the numerical results on the TIE from \autoref{subsec:ex2} we set $g=3\pi/2$ as in \autoref{eq:gbphi}. We also set $h=3\pi/2$. The viscosity parameter $\epsilon$ was chosen as $5e-2$. 
The stepsize in $z$-direction is $0.01$. 
The true phase, as given in \autoref{eq:gaussianbeam_seperated}, is
\[\phi(\vx,z) = -k\frac{\norm{\vx}^2}{2R(z)} + \operatorname{arg}\left( \frac{1}{q(z)}\right).\]
After transforming back to $\phi^{\epsilon}$, the numerical solution and error at $z = 0.1, 0.5, 1$ of the \autoref{eq:visc-bound} are shown in \autoref{fig:visco_rangez} and \autoref{fig:viscoerror_rangez}.

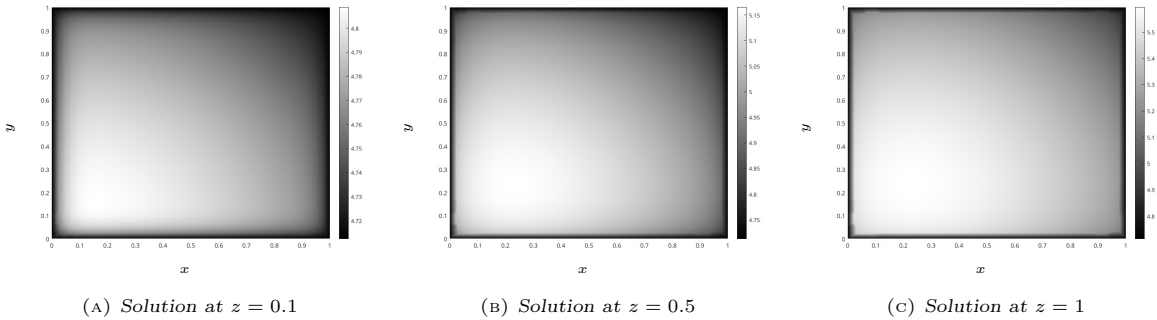
\begin{figure}[h]
    \centering
    \begin{minipage}{.32\textwidth}
    \begin{tikzpicture}
            \node (img) {\includegraphics[width=\textwidth]{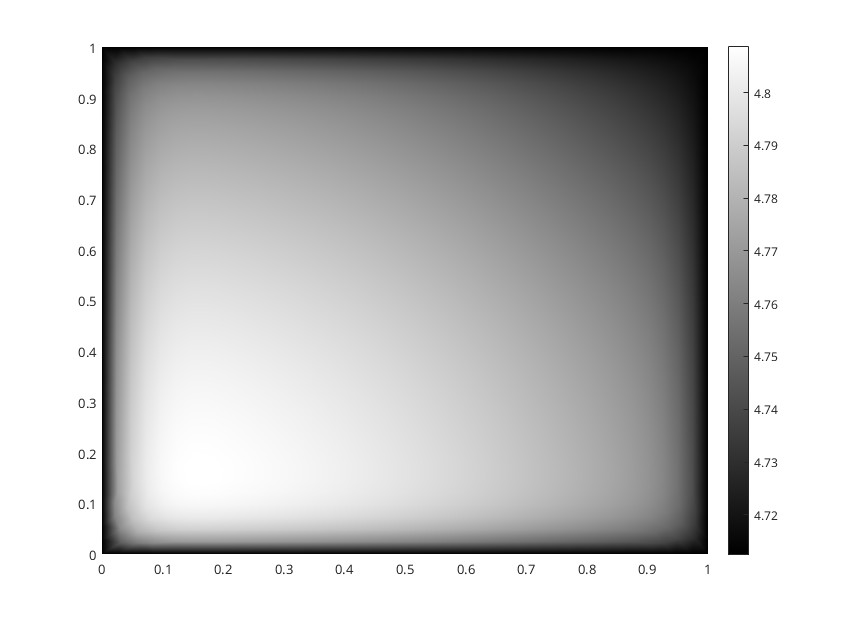}};
            \node[below=of img,yshift=1.3cm,xshift=-0.2cm]{\tiny$x$};
                                    \node[left=of img, rotate=90,anchor=center,yshift=-1.2cm]{\tiny $y$};
    \end{tikzpicture}
            \subcaption{Solution at $z = 0.1$}
            \label{fig:Solutionat0.1}
    \end{minipage}
    \begin{minipage}{.32\textwidth}
    \begin{tikzpicture}
            \node (img){ \includegraphics[width=\textwidth]{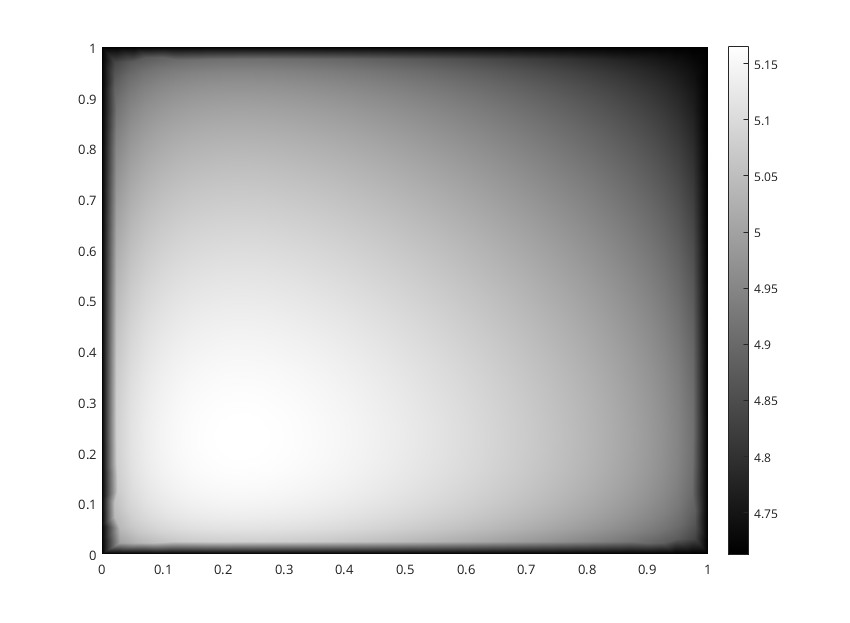}};
                        \node[below=of img,yshift=1.3cm,xshift=-0.2cm]{\tiny$x$};
                        \node[left=of img, rotate=90,anchor=center,yshift=-1.2cm]{\tiny $y$};
                \end{tikzpicture}
            \subcaption{Solution at $z= 0.5$}
            \label{fig:soluz.5}
    \end{minipage}
    \begin{minipage}{.32\textwidth}
    \begin{tikzpicture}
        \node (img) { \includegraphics[width=\textwidth]{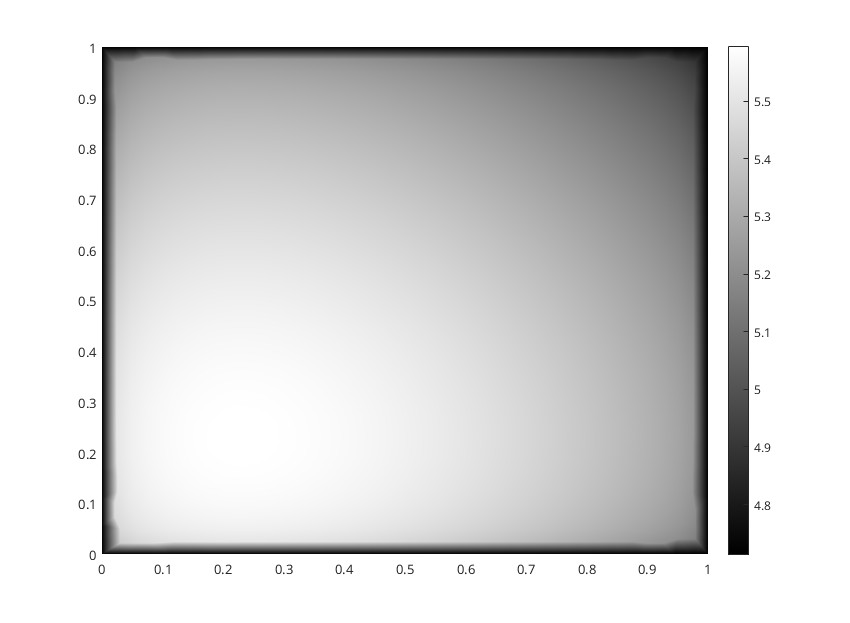}};
                    \node[below=of img,yshift=1.3cm,xshift=-0.2cm]{\tiny$x$};
                                            \node[left=of img, rotate=90,anchor=center,yshift=-1.2cm]{\tiny $y$};
    \end{tikzpicture}
         \subcaption{Solution at $z = 1$}
         \label{fig:solutionz1}
    \end{minipage}
    \caption{Snapshots of viscosity solutions for $\phi$ when $z = 0.1, 0.5, 1$ in \autoref{ex:viscosity_gauss}.  Given the result of the TIE in \autoref{fig:gb_correct_reconstruction} we used them as initial condition for the TPE. The results show the solution of the TPE in $z$-directions.}
    \label{fig:visco_rangez}
\end{figure}

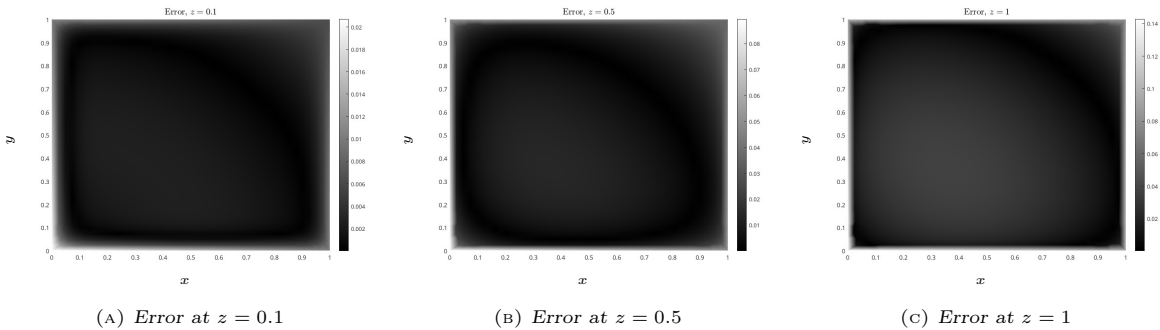
\begin{figure}[h]
    \centering
    \begin{minipage}{.32\textwidth}
    \begin{tikzpicture}
            \node (img) {\includegraphics[width=\textwidth]{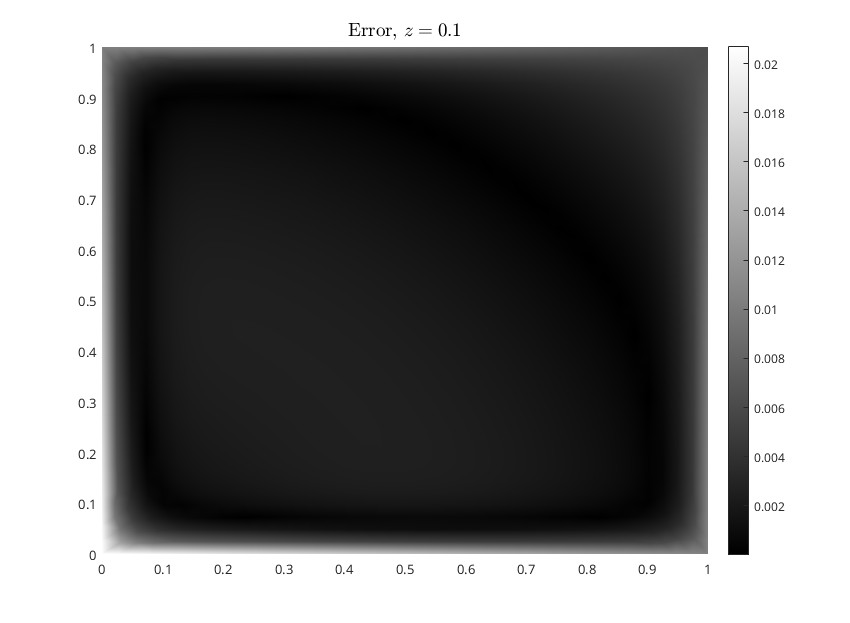}};
            \node[below=of img,yshift=1.3cm,xshift=-0.2cm]{\tiny$x$};
                                    \node[left=of img, rotate=90,anchor=center,yshift=-1.2cm]{\tiny $y$};
    \end{tikzpicture}
            \subcaption{Error at $z = 0.1$}
            \label{fig:Errorat0.1}
    \end{minipage}
    \begin{minipage}{.32\textwidth}
    \begin{tikzpicture}
            \node (img){ \includegraphics[width=\textwidth]{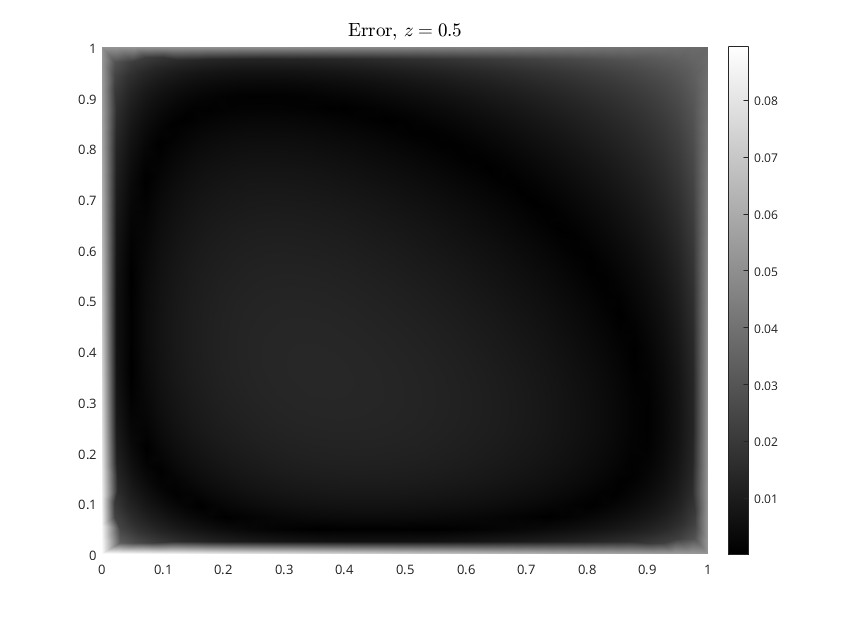}};
                        \node[below=of img,yshift=1.3cm,xshift=-0.2cm]{\tiny$x$};
                        \node[left=of img, rotate=90,anchor=center,yshift=-1.2cm]{\tiny $y$};
                \end{tikzpicture}
            \subcaption{Error at $z= 0.5$}
            \label{fig:errorz.5}
    \end{minipage}
    \begin{minipage}{.32\textwidth}
    \begin{tikzpicture}
        \node (img) { \includegraphics[width=\textwidth]{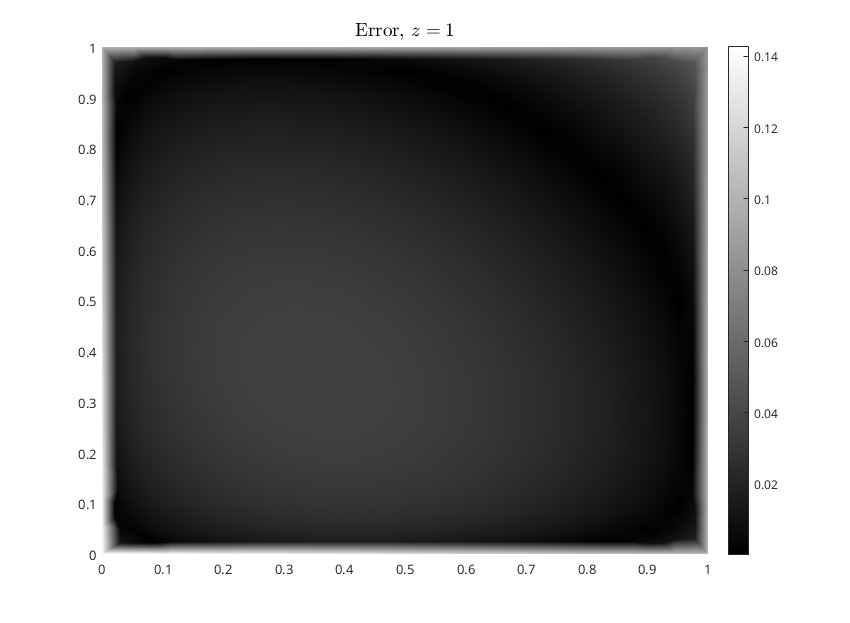}};
                    \node[below=of img,yshift=1.3cm,xshift=-0.2cm]{\tiny$x$};
                                            \node[left=of img, rotate=90,anchor=center,yshift=-1.2cm]{\tiny $y$};
    \end{tikzpicture}
         \subcaption{Error at $z = 1$}
         \label{fig:Errorz1}
    \end{minipage}
    \caption{Snapshots of absolute relative viscosity solution errors for $\phi$ when $z = 0.1, 0.5, 1$ in \autoref{ex:viscosity_gauss}, which increases approximately from $1\%$ to about $10\%$.}
    \label{fig:viscoerror_rangez}
\end{figure}
\end{example}
\begin{remark}
	We have investigated two methods for solving the nonlinear TPE, the characteristic method, which can be considered a theoretical tool, and the practical method of viscosity solutions. The latter, although computationally feasible, requires to specify virtual boundary conditions. These destroy potential conceptual benefits of the nonlinear TPE where solutions along characteristics can be computed without specifying boundary conditions (see \autoref{eq:colehopf}).   
\end{remark}

\section{Conclusion}
\label{sec:comparison}

In this paper, we derived the TIE and the TPE from the paraxial Helmholtz equation (see \autoref{eq:A}).
We considered two solutions of \autoref{eq:A} for testing the performance of the reconstructions.
The first example presented in \autoref{subsec:ex1} had constant intensity. 
As a second example we considered a Gaussian beam (see \autoref{subsec:ex2}).
We have reconstructed the phase given intensity measurements by solving the TIE in \autoref{sec:TIE}. For this we had to impose boundary conditions on $\Gamma_0$ and reconstruct the phase on the two-dimensional domain $\Omega_0$ as depicted in Case 1 of \autoref{fig:introboundary}.
It was shown that with ground-truth boundary conditions the reconstruction with the TIE works reasonably well, while imposing an arbitrary boundary condition resulted in significant deviations from the ground-truth. We presented the corresponding numerical results in \autoref{fig:phase:ex}, \autoref{fig:TIE:plane_wave}, \autoref{fig:TIE_numerical2} and \autoref{fig:TIE:gaussian}.

In \autoref{sec:TPE} we considered hybrid reconstruction of the phase in the 3-dimensional domain $\Omega$ following the strategy of Case 2 of \autoref{fig:introboundary}. Our current computational results are limited to the viscosity setting and the method of characteristic can be only used in exceptional, simple, cases.

Imposing virtual boundary condition is a computational necessity both for TIE and TPE, and might significantly affect the obtained reconstructions. Although theoretically the TPE has advantages in requiring no boundary conditions in exceptional cases, in computational practice they are still required as for the TPE. Moreover, after applying the Cole-Hopf transform the TPE is a second-order linear PDE, like the TIE.

\appendix
\subsection*{Acknowledgements}
This research was funded in whole, or in part, by the Austrian Science Fund
(FWF) 10.55776/P34981 -- New Inverse Problems of Super-Resolved Microscopy (NIPSUM) and
SFB 10.55776/F68 ``Tomography Across the Scales'', project F6807-N36
(Tomography with Uncertainties).  For the purpose of open access, the author has applied a CC BY public copyright license to any Author Accepted Manuscript version arising from this submission.
The financial support by the Austrian Federal Ministry for Digital and Economic
Affairs, the National Foundation for Research, Technology and Development and the Christian Doppler
Research Association is gratefully acknowledged.

\section*{References}
\renewcommand{\i}{\ii}
\printbibliography[heading=none]

\end{document}